\documentclass[12pt]{article}

\usepackage[usenames]{color}

\usepackage[all]{xy}
\def\xym{\xymatrix}

\def\chi{{\mathcal{X}}}

\def\rrd{{\rr^d}}

\def\cald{{\mathcal{D}}}
\def\calx{{\mathcal{X}}}

\def\({\left(}
\def\){\right)}

\def\vsp{\vspace*{1,5mm}\\ }

\def\bk{\bigskip }

\def\sk{\smallskip }

\def\n{\noindent }
\def\dd{\displaystyle}

\def\D{{\Delta}}

\def\barr{\begin{array}}
\def\earr{\end{array}}

\def\bit{\begin{itemize}}

\def\eit{\end{itemize}}

\def\D{{\Delta}}
\def\FP{Fokker--Planck}

\usepackage{amsmath, amsthm, amscd, amsfonts, amssymb}

\newtheorem{theorem}{Theorem}[section]
\newtheorem{proposition}[theorem]{Proposition}

\newtheorem{claim}[theorem]{Claim}

\theoremstyle{definition}

\newtheorem{remark}[theorem]{Remark}

\def\1{^{-1}}
\def\vsp{\vspace*{2mm}\\ }

\def\calf{{\mathcal{F}}}

\def\calx{{\mathcal{X}}}

\def\rr{{\mathbb{R}}}

\def\9{{\infty}}

\def\a{{\alpha}}
\def\b{{\beta}}

\def\g{{\gamma}}
\def\wt{\widetilde}

\def\vf{{\varphi}}

\def\pp{{\partial}}
\def\D{{\Delta}}
\def\vp{{\varepsilon}}
\def\barr{\begin{array}}
\def\earr{\end{array}}

\def\dd{\displaystyle}
\def\bk{\bigskip }
\def\sk{\smallskip}
\def\n{\noindent }

\def\pas{\mathbb{P}\mbox{-a.s.}}

\def\vsp{\vspace*{2mm}\\ }
\def\ff{\forall }
\def\({\left(}
\def\){\right)}
\def\<{\left<}
\def\>{\right>}

\title{Uniqueness for nonlinear Fokker--Planck equations and weak uniqueness for~McKean-Vlasov SDEs} 
\author{Viorel Barbu\thanks{Octav Mayer Institute of Mathematics of the Romanian Academy,     Ia\c si, Romania.  Email: vbarbu41@gmail.com}\and Michael R\"ockner\thanks{Fakult\"at f\"ur Mathematik, Universit\"at Bielefeld,  D-33501 Bielefeld, Germany.  Email: roeckner@math.uni-bielefeld.de}\ \thanks{Academy of Mathematic and System Sciences, CAS, Beijing, China}}
\date{}

\begin{document}
\maketitle
\begin{abstract}
\n One proves the uniqueness of distributional solutions to  nonlinear Fokker--Planck equations with monotone diffusion term and derive as a  consequence (restricted) uniqueness in law for the correspon\-ding  McKean--Vlasov stochastic differential equation (SDE).\sk\\
{\bf Mathematics Subject Classification (2000):} 60H30, 60H10, 60G46, 35C99.\\
{\bf Keywords:} Fokker--Planck  equation, mild solution, distributional solution.   
\end{abstract}

\section{Introduction}\label{s1}

Consider  the nonlinear \FP\ equation  
\begin{equation}\label{e1.1}
\barr{ll}
u_t-\Delta\beta(u)+{\rm div}(b(x,u)u)=0\mbox{ in }\cald'((0,\9)\times\rr^d),\vsp 
u(0,x)=u_0(x),
\earr\end{equation}
where $\b:\rr\to\rr$ and $b:\rrd\times\rr\to\rrd$ satisfy  the following assumptions
\begin{itemize}
	\item[(i)] {\it$\beta(0)=0$, $\b\in C^1(\rr)$, and 
		\begin{equation}
		\g_0|r_1-r_2|^2\le(\b(r_1)-\b(r_2))(r_1-r_2),\ r_1,r_2\in\rr,\label{e1.2}
			\end{equation}
	where $0<\g_0<\9.$	}
	\item[(ii)] $b\in C_b(\mathbb R^{d+1}; \rrd)$, $b(x, \cdot) \in C^1(\rrd; \rrd)$ $ \forall x\in\rrd$ {\it such that}  $$\mbox{$\sup\{|b^i_r(x,r)|;x\in\rrd,i=1,2,$ $|r|\le M\}\le C_M,$ $\ff M>0.$}$$
	\item[(ii')] $b(x,0)=0$ $\forall x\in\rrd$, $b \in C^1(\mathbb R^{d+1}; \rrd) $ {\it and for
		$$\delta(r):=\sup\{|b_x(x,r)|;\ x\in\rrd\},$$we have} $\delta\in   C_b(\rr).$
\end{itemize}
Here$$b(x,u)=\{b^i(x,u)\}^d_{i=1}\mbox{ and }b^i_r=\frac{\pp b^i}{\pp r}, b_x=\big\{ \nabla_x \, b^i(x, \cdot) \big\}_{i=1}^d.$$
By a distributional solution (in the sense of Schwartz) with initial condition $u_0\in L^1$ we mean a function $u:[0,\9)\to L^1(\rrd)$ such that
$(u(t,\cdot)dx)_{t\in[0,T]}$ is narrowly continuous, that is,
\begin{eqnarray}
&&\lim_{t\to s}\int_{\rrd}u(t,x)\psi(x)dx=\int_{\rrd}u(s,x)\psi(x)dx,\ \ff\psi\in C_b(\rrd),\ s\ge0,\qquad\label{e1.3}\\[2mm]
&&\dd\int^\9_0\int_{\rr^d}(u(t,x)\vf_t(t,x)+\b(u(t,x))\D\vf(t,x)\label{e1.4}\\
&&\qquad+b(x,u(t,x))u(t,x))\cdot\nabla_x\vf(t,x))dt\,dx=0,\nonumber\\
&&\qquad\qquad\qquad\qquad\qquad\qquad\ff\vf\in C^\9_0((0,\9)\times\rr^d).\nonumber
\end{eqnarray}
(In the following, we shall use the notation $b^*(x,u)=b(x,u)u.$)

In \cite{1} it was proved, in particular, that, if (i)--(ii) hold and, in addition, for $\Phi(u)\equiv\frac{\beta(u)}u,$ $u\in\rr,$ we have $\Phi\in C^2(\rr)$, then there is a  mild solution $u\in C([0,\9);L^1(\rr^d))$ for each $u_0\in L^1(\rr^d)$. The mild solution $u$ is defined~as 
$$u(t)=\lim_{h\to0}u_h(t)\mbox{ in }L^1(\rr^d),\ \ff t\ge0,$$where
\begin{equation}
\label{e1.6}
\barr{l}
u_h(t)=u^i_h\mbox{ for }t\in[ih,(i+h)h],\ i=0,1,...,Nh=T,\vsp 
u^{i+1}_h-h\D\b(u^{i+1}_h)+h\,{\rm div}(b(x,u^{i+1}_h)u^{i+1}_h)=u^i_h\mbox{ in }\cald'(\rr^d),\\\hfill i=0,1,...,\vsp 
u^0_h=u_0.\earr	
\end{equation}Moreover, $S(t)u_0=u(t),\ t\ge0,$ is a strongly continuous semigroup of nonexpansive mappings in $L^1(\rrd)$. 

As easily seen, any mild solution is a distributional solution but the uniqueness  follows in the class of mild solutions only. Here, we shall prove the uniqueness for \eqref{e1.1} in the class of distributional solutions and derive from  this result the uniqueness in law of solutions to McKeen--Vlasov SDE
\begin{equation}
\label{e1.7}
dX(t)=b(X(t),u(t,X(t)))dt+\frac1{\sqrt{2}}\(\frac{\b(u(t,X(t))}{u(t,X(t))}\)^{\frac12}dW(t).
\end{equation}

\n{\bf Notation.}  Denote by $L^p(\rr^d)=L^p$ the  space of $p$-summable functions on $L^p$, with the norm denoted $|\cdot|_p$. By $H^k(\rrd)=H^k$, $k=1,2$, and $H^{-k}(\rrd)=H^{-k},$ we denote the standard Sobolev spaces on $\rrd$ and by $C_b(\rrd)$ the space of continuous and bounded functions on $\rrd$. By $C^k(\rrd)$ we denote the space of continuously differentiable functions on $\rrd$ of order $k$, by $C^1_b(\rrd)$ the space $\left\{u\in C^1(\rrd);\frac{\pp u}{\pp y}\in C_b(\rrd),\ j=1,...,d\right\}.$ The spaces of continuous and differentiable functions on $(0,T)\times\rrd$ are denoted in a similar way and we shall simply write
$$C^1_b(\rrd)=C^1_b,\ C^k(\rrd)=C^k,\ k=1,2.$$ The scalar product in $L^2$ is denoted $\<\cdot,\cdot\>_2$ and by ${}_{H^{-1}}\!\<\cdot,\cdot\>_{H^1}$ the pairing between $H^1$ and $H^{-1}.$ Of course, on $L^2\times L^2$ this coincides with $\<\cdot,\cdot\>_2.$ The~scalar product $\<\cdot,\cdot\>_{-1}$ on $H\1$ is  taken~as
\begin{equation}
\label{e1.8}
\<u,v\>_{-1}=\<(I-\D)\1u,v\>_2,\ \ff u,v\in H\1
\end{equation}with the corresponding norm
\begin{equation}
\label{e1.9}
|u|_{-1}=(\<u,u\>_{-1})^{\frac12},\ u\in H\1.
\end{equation}
By $\cald'((0,\9)\times\rrd)$ and $\cald'(\rrd)$ we denote the space of Schwartz distributions on $(0,\9)\times\rrd$ and $\rrd$, respectively. If $\calx$ is a Banach space, we denote by $W^{1,2}([0,T];\calx)$ the infinite dimensional Sobolev space 
$\{y\in L^2(0,T;\calx);$ $\frac{dy}{dt}\in L^2(0,T;\calx)\}$, where $\frac d{dt}$ is taken in the sense of vectorial distributions. We also set, for each $z\in C^1(\rrd\times\rr)$,
$$z_r(x,r)=\frac\pp{\pp r}\ z(x,r),\ \ z_x=\nabla_x z(x,r).$$
We shall denote the norms on $\rrd$ and $\rr$  by the same symbol $|\cdot|$.

\section{The main result}
\setcounter{equation}{0}

The next result is a uniqueness theorem for distributional solutions $u$ to \eqref{e1.1}. In the special case $b\equiv0$, such a uniqueness result for\eqref{e1.1}  was established earlier in \cite{3} for continuous and monotonically nondecreasing functions $\b$.  

\begin{theorem}\label{t1} Let $T>0$ and let conditions {\rm(i)--(ii)} on $\beta$ and $b$ hold. For each $u_0\in L^\9\cap L^1$, the \FP\ equation \eqref{e1.1} has at most one distributional solution 
\begin{equation}\label{add-e5}
u\in L^1((0,T);L^1)\cap L^\9((0,T)\times\rrd).
\end{equation}
\end{theorem}

\n{\bf Proof.} Let $u_1,u_2\in L^1(0,T;L^1)\cap L^\9((0,T)\times\rrd)$ be two distributional solutions to \eqref{e1.1} and let $u=u_1-u_2.$ We have
\begin{equation}
\label{e2.1}
\barr{l}
u_t-\D(\b(u_1)-\b(u_2))+{\rm div}(b^*(x,u_1)-b^*(x,u_2))=0\\\hfill\mbox{ in }\cald'((0,\9)\times\rrd)\vsp 
u(0,x)=0.\earr
\end{equation}
(Here, $b^*(x,r)=b(x,r)r,\ \ff x\in\rrd,\ r\in\rr.$)

It should be mentioned that, by (i), (ii), it follows that $u_i,\beta(u_i), b^*(\cdot, u_i)\in L^2((0,T);L^2),$ $i=1,2$, and, therefore, $u\in W^{1,2}([0,T];H^{-2}).$ 

Consider the operator $\Gamma:H\1\to H^1$ defined by
$$\Gamma u=(1-\D)^{-1}u,\ u\in H\1(\rrd)$$and note that $\Gamma$ is an isomorphism of $H\1$ onto $H^1$ and $\Gamma\in L(H^{-2},L^2).$
 Since $u_i\in L^2(0,T;L^2),$ $i=1,2,$ it follows that $y=\Gamma u\in L^2(0,T;H^2)\cap W^{1,2}([0,T];L^2)$ and so, by~\eqref{e2.1},  we have
\begin{equation}
\label{e2.2}
\barr{l}
\dd\frac{dy}{dt}-\Gamma\D(\b(u_1)-\b(u_2))+\Gamma\,{\rm div}(b^*(x,u_1)-b^*(x,u_2))=0,\\\hfill\mbox{ a.e. }t\in(0,T),\\
y(0)=0,\earr
\end{equation}
where $\frac{dy}{dt}\in L^2(0,T;L^2).$  (We note that here $\frac{dy}{dt}$ is taken in the sense of $L^2$-valued vectorial distributions on $(0,T)$ and so $y:[0,T]\to L^2$ is absolutely continuous.)

Assume first that, in addition,
\begin{equation}
\label{add-e6}
u\in L^2(0,T;H^1);\ u\in W^{1,2}([0,T];H^{-1}).
\end{equation}

Now, we take the scalar product in $L^2$ of \eqref{e2.2}  with $u=u_1-u_2$. Taking into account that
$$\<\frac{du}{dt}\,(t),u(t)\>_{-1}=\frac12\ \frac d{dt}\,|u(t)|^2_{-1}, \mbox{ a.e. }t\in(0,T),$$we get, by \eqref{e2.2} that
$$\barr{l}
\dd\frac 12\ \frac d{dt}\,|u(t)|^2_{-1}+\<\b(u_1)-\b(u_2),u_1-u_2\>_2
=\<\Gamma(\b(u_1)-\b(u_2)),u_1-u_2\>_2\vsp 
\qquad-\<\Gamma\,{\rm div}((b^*(x,u_1)-b^*(x,u_2)),u_1-u_2\>_2,\ \mbox{a.e. }t\in(0,T).\earr$$ 
By \eqref{e1.2}, this yields for a.e.~$t\in(0,T)$
\begin{equation}
 \label{e2.3}
\barr{r}
\dd\frac12\ \frac d{dt}\,|u(t)|^2_{-1}+\g_0|u(t)|^2_2
\le\<\b(u_1(t))-\b(u_2(t)),u_1(t)-u_2(t)\>_{-1}\vsp
\qquad-\<{\rm div}(b^*(x,u_1(t))-b^*(x,u_2(t))),u_1-u_2\>_{-1}.\earr
\end{equation}
We note that
\begin{equation}
 \label{e2.4}
 |\Gamma f|_2\le|f|_2,\ \ff f\in L^2.
\end{equation}
We also have
\begin{equation}
 \label{e2.5}
|{\rm div}\,F|_{-1}\le2|F|_2,\ \ff F\in (L^2)^d.
\end{equation}This  yields
\begin{equation}
 \label{e2.6}
 \barr{ll}
|\<\b(u_1)-\b(u_2),u_1-u_2\>_{-1}| \le|\b(u_1)-\b(u_2)|_{2}|u|_{-1} 
\le\b_M|u|_2|u|_{-1}
\earr\end{equation}and
\begin{equation}
 \label{e2.7}
\barr{l}
\left|\<{\rm div}\ (b^*(x,u_1)-b^*(x,u_2)),u_1-u_2\>_{-1}\right|\vsp
\qquad\le
2|b^*(x,u_1)-b^*(x,u_2)|_2|u_1-u_2|_{-1}\vsp
\qquad\le2(|b|_\9+b_M|u_1|_\9)|u|_2|u|_{-1},
\earr
\end{equation}
where $M=\max\{|u_1|_\9,|u_2|_\9)$ and
$$\barr{lcl}
\b_M&=&\dd\sup\{\b'(r);\ |r|\le M\},\vsp 
b_M&=&\dd\sup\left\{\frac{|b(x,u_1)-b(x,u_2)|}{|r_1-r_2|};\ x\in\rrd,\ |r_1|,|r_2|<M\right\}\vsp&\le&
\sup\{|b_r(x,r)|;|r|\le M,\ x\in\rrd\} < \infty \earr$$
by (ii).

By \eqref{e2.3}--\eqref{e2.7}, we get
$$\barr{r}
\dd\frac12\ \frac d{dt}\,|u(t)|^2_{-1}+\g_0|u(t)|^2_{2}
\le(\b_M+2(|b|_\9+b_M|u_1|_\9)
|u(t)|^{}_2|u(t)|^{}_{-1},\\\hfill\mbox{a.e. }t\in(0,T).\earr$$
This yields
$$\frac d{dt}\,|u(t)|^2_{-1}\le C|u(t)|^2_{-1},\mbox{ a.e. }t\in(0,T).$$Since $u:[0,T]\to H^{-1}$ is absolutely continuous and narrowly continuous, we infer that $|u(t)|_{-1}=0$, $\ff t\in[0,T]$, and so $u\equiv0$, as claimed.
To conclude the proof, we shall prove

\begin{claim}\label{add-l1} If $u$ is a distributional solution to \eqref{e1.1} satisfying \eqref{add-e5}, then \eqref{add-e6} holds.
\end{claim}

\n{\bf Proof of Claim \ref{add-l1}.} $1^\circ$. We shall first assume   that $u$ satisfies   the stronger initial condition
\begin{equation}
\label{add-e7}
\mbox{essential }\lim_{t\to0} \int_{\rr^d} |u(t,x)-u_0(x)|dx=0.
\end{equation}
Next, we consider the truncated function
\begin{equation}
\label{add-e8}
\beta_M(r)=\left\{\barr{ll}
\beta(r)&\mbox{ if }|r|\le M,\vsp
\beta(M)+\beta'(M)(r-M)&\mbox{ if }r>M,\vsp
\beta(-M)+\beta'(-M)(r+M)&\mbox{ if }r<-M,\earr\right.
\end{equation}where $M=\|u\|_{L^\9(0,T)\times\rr^d)}.$

We set
$$f(t,x)\equiv b^*(x,u(t,x)),\ \ f_\nu=f*\rho_\nu,\ \ \nu>0,$$where $\rho_\nu$ is a standard mollifier on $\rr^{d+1}$ and $b^*(x,u)\equiv b(x,u)u$.

We note that in equation \eqref{e1.1} one can take $\beta\equiv\beta_M$ because $|u|_\9\le M.$

We consider the equations
\begin{equation}
\label{add-e9}
\barr{l}
\wt u_t-\Delta\beta_M(\wt u)+{\rm div}(f)=0\mbox{ in }(0,T)\times\rr^d,\vsp
\wt u(0)=u_0,\earr\qquad
\end{equation}
\begin{equation}
\label{add-e10}
\barr{l}
(\wt u_\nu)_t-\Delta\beta_M(\wt u_\nu)+{\rm div}(f_\nu)=0\mbox{ in }(0,T)\times\rr^d,\vsp
\wt u_\nu(0)=u_0.\earr
\end{equation}
Clearly, \eqref{add-e9} has a unique solution $\wt u$ satisfying \eqref{add-e6}, i.e.,
\begin{equation}
\label{add-e11}
\wt u \in L^2(0,T;H^1),\ \wt u\in W^{1,2}([0,T];H^{-1}),
\end{equation}
while $\wt u_\nu$ is more regular, that is,   additionally to \eqref{add-e11} it satisfies
\begin{equation}\label{add-e12}
\wt u_\nu\in L^\9((0,T)\times \rr^d)\cap C^1([0,T];L^1),
\end{equation}because the operator $u\to-\Delta\beta_M(u)$ is $m$-accretive in $L^1(\rr^d)$ and $u_0\in(L^1\cap L^\9)(\rr^d)$, ${\rm div}\,f_\nu\in(L^1\cap L^\9)((0,T)\times\rr^d)$. 
It is also clear that
\begin{equation}
\label{add-e13}
\barr{ll}
\dd\lim_{\nu\to0} \wt u_\nu=\wt u
&\mbox{ strongly in }L^2(0,T;H\1).\earr
\end{equation}
To prove the claim, we are going to show that
\begin{equation}
\label{add-e14}
\wt u=u,\mbox{ a.e. in }(0,T)\times\rr^d.
\end{equation}
To this end, we shall invoke an argument due to Brezis \& Crandall \cite{3}.

 Namely, we subtract equations \eqref{e1.1}, \eqref{add-e10} and get
 $$(\wt u_\nu- u)_t-\Delta(\beta(\wt u_\nu)-\beta(u)=-{\rm div}(f_\nu-f)\mbox{ in }\cald'((0,T)\times\rr^d),$$
 $$\mbox{essential }\lim_{t\to0}\int_{\rr^d}|\wt u_\nu(t,x)-u(t,x)|dx=0.$$
 We set $z_\nu=\wt u_\nu-u$, $h_\nu=\beta(\wt u_\nu)-\beta(u)$ and get
 \begin{equation}
 \label{add-e14a}
 (z_\nu)_t-\Delta h_\nu=-{\rm div}(f_\nu-f)\mbox{   in }\cald'((0,T)\times\rr^d),
 \end{equation}
\begin{equation}
\label{add-e14'}
\mbox{essential }\lim_{t\to0}\int_{\rr^d}|z_\nu(t,x)|dx=0.
\end{equation}
Arguing as in \cite{3} (proof of Proposition 1), we set
$$g^\nu_\vp(t)=(B_\vp(z_\nu),z_\nu),\ B_\vp=(\vp I-\Delta)^{-1},$$where $(\ ,\ )$ denotes the inner product in $L^2=L^2(\rr^d)$, and get
$$(g^\nu_\vp(t))_t=2(\vp B_\vp(h_\nu)-h_\nu,z_\nu(t))+2(f_\nu-f,\nabla B_\vp(z_\nu))$$and, therefore (we note that $(h_\nu,z_\nu)\ge0$),
 \begin{equation}
\label{add-e15}
g^\nu_\vp(t)\le2\dd\int^t_0\vp(B_\vp(h_\nu(s)),z_\nu(s))ds
+2\dd\int^t_0((f_\nu-f)(s),\nabla B_\vp(z_\nu(s)))ds,
\end{equation}because, by \eqref{add-e14'}, we have, since $z_\nu\in L^\9\cap L^1$,
 \begin{equation}
\label{add-e16}
g^\nu_\vp(0^+)=\mbox{essential }\lim_{t\to0}(B_\vp(z_\nu(t),z_\nu(t))=0.
\end{equation}
Recalling that, for a.e. $s\in[0,T]$,
$$\vp B_\vp(z_\nu(s))-\Delta B_\vp(z_\nu(s))=z_\nu(s),$$we get
$$g^\nu_\vp(s)=(B_\vp(z_\nu(s)),z_\nu(s))=\vp|B_\vp(z_\nu(s))|^2_{L^2}+|\nabla B_\vp(z_\nu(s))|^2_{L^2}$$
and so, by \eqref{add-e15}, we have
 \begin{equation}
\label{add-e17}
g^\nu_\vp(t)\le2\dd\int^t_0\vp(B_\vp h_\nu(s),z_\nu(s))ds
+\dd\int^t_0|(f_\nu-f)(s)|^2_{L^2}ds+\dd\int^t_0 g^\nu_\vp(s)ds.
\end{equation}Hence, by Gronwall's inequality,
 \begin{equation}
\label{add-e18}
g^\nu_\vp(t)\le\int^t_0(2\vp(B_\vp (h_\nu(s)),z_\nu(s))+|(f_\nu-f)(s)|^2_{L^2})e^{t-s}ds.
\end{equation}
On the other hand, by Lemma 1 in \cite{3} we have for a.e. $s\in[0,T]$ and also in $L^1(0,T)$
 \begin{equation}
\label{add-e18a}
\lim_{\vp\to0}(\vp B_\vp (h_\nu(s)),z_\nu(s))=0.
\end{equation}
By \eqref{add-e18}, this yields
 \begin{equation}
	\label{add-e22a}
	(G*z_\nu(t),z_\nu(t))\le\int^t_0|(f_\nu-f)(s)|^2_{L^2}e^{t-s}ds,\end{equation}
since $$\lim_{\vp\to0}g^\nu_\vp(t)=((-\Delta)^{-1}z_\nu(t),z_\nu(t))=(G*z_\nu(t),z_\nu(t)),$$where $G$ is the classical Newtonian Green function. We note that the  latter  inner product is well defined, since $z_\nu\in\dd\bigcap_{p\in[1,\9]}L^p$ and hence $G*z_\nu\in L^{\frac{2d}{d-2}}.$ Letting $\nu\to0$ in \eqref{add-e22a}, we find
$$\lim_{\nu\to0}(G*z_\nu(t),z_\nu(t))=0.$$
In particular, $\xym{z_\nu(t)\ar[r]_{\ \ \nu\to0}&0}$ in $H^{-1}$. Since $\xym{z_\nu\ar[r]_{\nu\to0\ \ }&\wt u-u}$ in $L^2(0,T;L^2)$ by~\eqref{add-e13}, we have along a subsequence $\xym{z_\nu(t)\ar[r]_{\nu\to0\ \ \ \ }&(\wt u-u)(t)}$ in $L^2$ for a.e. $t\in[0,T]$, hence \eqref{add-e14} follows.

\bk\n$2^\circ$. {\it The proof without assuming  \eqref{add-e7}.} It follows by the argument described under (1.12)  in \cite{3} that \eqref{e1.4} indeed implies \eqref{add-e16}. So, assumption \eqref{add-e7} can be dropped.

This completes the proof.
\hfill$\Box$\bk 

\n As mentioned earlier, under Hypotheses (i)--(ii), if $\Phi\in C^2$, where $\Phi(u)\equiv\frac{\beta(u)}u$, $u\in\rr$, then the \FP\ equation \eqref{e1.1}, for each $u_0\in L^1$, has a unique mild solution $u\in C([0,\9);L^1).$ This mild solution is also easily checked to be a distributional solution to \eqref{e1.1}. As regards this solution, we also have

\begin {proposition}\label{p2} Assume that {\rm (i), (ii), (ii')} hold, and that, for $\Phi(u)\equiv\frac{\beta(u)}u$,
\begin{itemize}
	\item[{\rm(iii)}] $\Phi\in C^2(\rrd).$
	\end{itemize}
	 Then, for each $u_0\in L^1\cap L^\9$, the mild solution $u$ to \eqref{e1.1} satisfies also
\begin{equation}
 \label{e2.8}
u\in L^\9((0,T)\times\rrd),\ \ff T>0.
\end{equation}
\end{proposition}

\n{\bf Proof.} We rewrite \eqref{e1.1} as
\begin{equation}
 \label{e2.9}
\barr{l}
(u-|u_0|_\9-\alpha(t))_t-\D(\b(u)-\b(|u_0|_\9+\a(t)))\vsp \qquad
+{\rm div}(b^*(x,u)-b^*(x,|u_0|_\9+\a(t)))\vsp 
\qquad=-{\rm div}(b^*(x,|u_0|_\9+\a(t)))-\a'(t)\le0\mbox{ in }(0,\9)\times\rrd,
\earr
\end{equation}where $\a\in C^1([0,\9))$ is chosen in such a way that
\begin{equation}
\label{e2.10}
\hspace*{-3mm}\barr{l}
\a'(t)+\sup\{|b_x(x,|u_0|_\9+\a(t))|;x\in\rrd\}(|u_0|_\9+\a(t))=0,\   t\in(0,T),\vsp 
\a(0)=0.\earr
\end{equation}
We may find $\a$ of the form $\a=\eta-|u_0|_\9,$ where $\eta$ is a solution to the equation 
\begin{equation}
\label{e2.11}
\barr{l}
\eta'-\delta(\eta)\eta  =0,\ t\ge0,\vsp 
\eta(0)=|u_0|_\9,\earr
\end{equation} $\delta(r)=\sup\{|b_x(x,r)|;x\in\rrd\},$ $r\in\rr.$ Clearly, \eqref{e2.11} has such a solution    $\eta\in C^1([0,\9)),$ $\eta\ge0,$ on $[0,\9)$ because $\delta\in C_b(\rr).$

Formally, if we multiply \eqref{e2.9} by ${\rm sign}(u-|u_0|_\9-\a)^+$,   integrate over $\rrd$  and use the monotonicity of $\b$, we get by \eqref{e2.5} that  
\begin{equation}\label{e2.11a}
	\frac d{dt}\,|(u(t)-|u_0|_\9-\a(t))^+|_1\le0,\ \mbox{a.e. }t\in(0,T). \end{equation}This yields $u(t)\le|u_0|_\9+\a(t)$, $\ff t\ge0$, and similarly it follows that $u(t)\ge-|u_0|_\9-\a(t).$ Hence, $u\in L^\9((0,T)\times\rrd),$ as claimed.

The above formal argument can be made rigorous if $u$ is a strong solution to \eqref{e1.1} (which is not the case here). Then (see the detailed argument in \cite{1})
\begin{equation}\label{e2.11aa}
\hspace*{-2mm}\barr{l}
\dd\lim_{\delta\to0}\frac1\delta
\int_{[0<(\beta(u)-\beta(|u_0|_\9+\a(t))^+)\le\delta]}
|b^*(x,u)-b^*(x,|u_0|_\9+\a(t))|\,|\nabla u|dx\vsp
\qquad=\dd\lim_{\delta\to0}\frac1\delta
\int_{[0<(\beta(u)-\beta(|u_0|_\9+\a(t))^+)\le\delta]}
(|b(x,u)-b(x,|u_0|_\9+\a(t))|\,|u|
\vsp\qquad+|b(x,|u_0|_\9+\a(t))|\,|u-|u_0|_\9-\a(t)|)|\nabla u|dx=0,\ \ff t\in(0,T),\earr\end{equation}which is true if $\nabla u\in L^2(0,T;L^2)$ and $b(x,\cdot)\in{\rm Lip}(\rr)$ uniformly in $x$ (which is the case if $b_u\in C_b(\rrd\times\rr))$. In order to be in  such a situation, we approximate \eqref{e1.1}~by
\begin{equation}\label{e2.11aaa}
\barr{l}
u_t-\D(\b(u)+\vp\b(u)+{\rm div}(b_\vp(x,u)u))=0\mbox{ in }(0,T)\times\rrd,\vsp
u(0,x)=u_0(x),\earr\end{equation}where $\vp>0$ and $b_\vp\in C^1_b(\rrd\times\rr)$ is a smooth approximation of $b$. (For instance, $b_\vp=b*\rho_\vp$, where $\rho_\vp$ is a standard mollifier.) Then, as proved earlier in \cite{1}, \cite{2az}, \cite{2}, equation \eqref{e2.11aaa} has a unique solution $u_\vp\in L^2(0,T;H^1)\cap C([0,T];L^1)\cap W^{1,2}([0,T];H\1)$ and $u_\vp\to u$ in $C([0,T];L^1)$ as $\vp\to0$. An easy way to prove this is to apply the Trotter--Kato theorem to the family of $m$-accretive operators in $L^1$
$$\barr{lcl}
A_\vp u&=&-\D\b(u)+\vp\b(u)+{\rm div}(b_\vp(x,u)u),\vsp 
D(A_\vp)&=&\{u\in L^1;-\Delta\beta(u)+\vp\beta(u)+{\rm div}(b_\vp(x,u)u)\in L^1\}.\earr$$
(See the argument in \cite{2}.) Then, we replace \eqref{e2.9} by
\begin{equation}\label{e2.14a}
\barr{l}
(u_\vp-|u_0|_\9-\a(t))_t-\D(\b(u_\vp)-\b(|u_0|_\9+\a(t)))\vsp
+\vp(\b(u)-\b(|u_0|_\9+\a(t)))+{\rm div}(b^*_\vp(x,u_\vp)-b^*_\vp(x,|u_0|_\9+\a(t))\vsp
=-b^*_\vp(x,|u_0|_\9+\a(t))-\a'(t)-\vp\beta(|u_0|_\9+\a(t))\le0,\vsp\hfill \mbox{ a.e. in }(0,T)\times\rrd,\earr\end{equation}
where $b^*_\vp(u)=b_\vp(u)u.$

Let $\calx_\delta\in{\rm Lip}(\rr)$ be the following approximation of the signum function
$$\calx_\delta(r)=\left\{\barr{rll}
1&\mbox{for }&r\ge\delta,\vsp 
\dd\frac r\delta&\mbox{for }&|r|<\delta,\vsp 
-1&\mbox{for }r<-\delta,\earr\right.$$where $\delta>0.$ If we multiply \eqref{e2.14a} by $\calx_\delta((\beta(u_\vp)-\beta(|u_0|_\9+\a))^+)$ and integrate over $\rrd$, we get
$$\!\!\!\barr{l}
\dd\int_{\rrd}(u_\vp-|u_0|_\9-\a)_t\calx_\delta((\beta(u_\vp)-\beta(|u_0|_\9+\a))^+)dx\\ 
\qquad\le\dd\frac1\delta
\int_{[0<(\beta(u_\vp)-\beta(|u_0|_\9+\a))^+\le\delta]}
(b^*(x,u_\vp)u_\vp-b^*_\vp(x,|u_0|_\9+\a))\cdot\nabla u_\vp\,dx,\\\hfill \ff t\in(0,T),\earr$$because $\beta$ is monotonically increasing and 
$$\barr{r}
\dd\nabla(\beta(u_\vp)-\beta(|u_0|_\9+\a)\cdot\nabla\calx_\delta
((\beta(u_\vp)-\beta(|u_0|_\9+\a))^+)\ge0\mbox{ in } (0,T)\times\rrd.\earr$$
 Then, by \eqref{e2.11aa}, we get, for $\delta\to0$,
$$\int_\rrd(u_\vp-|u_0|_\9-\a(t))^+_t\,dx\le0,\ \ff t\in(0,T),$$and this yields
$$u_\vp(t,x)-|u_0|_\9-\a(t)\le0,\mbox{ a.e. on }(0,T)\times\rrd,$$and so, $u_\vp\le|u_0|_\9+\a$, a.e. on $(0,T)\times\rrd$. Then, we pass to the limit $\vp\to0$ to get the claimed inequality.\hfill$\Box$\bk

By Theorem \ref{t1} and Proposition \ref{p2}, we therefore get the following existence and uniqueness result for \eqref{e1.1}.  

\begin{theorem}\label{t4} Under hypotheses {\rm(i), (ii), (ii'), (iii)}, for each $u_0\in L^1\cap L^\9$, equation \eqref{e1.1} has a unique distributional solution
\begin{equation}
\label{e2.12}
u\in L^1((0,T);L^1)\cap L^\9((0,T)\times\rrd),\ \ff T>0.
\end{equation}
\end{theorem}

\section{Uniqueness of the linearized equation}
\setcounter{equation}{0}

Consider a distributional solution of the linearized equation corresponding to \eqref{e1.1}, that is,
\begin{equation}
\label{e3.1}
\barr{l}
v_t-\D(\Phi(u)v+{\rm div}(b(x,u)v)=0\mbox{ in }\cald'((0,\9)\times\rrd,\vsp 
v(0,x)=v_0(x),\earr
\end{equation}
where $u\in L^\9((0,T)\times\rrd),\ \ff T>0$. By (i)--(ii), we have
$$b(x,u),\Phi(u)=\frac{\b(u)}{u}\in L^\9((0,\9)\times\rrd).$$Moreover, we have
\begin{equation}
\label{e3.2}
\Phi(u)\ge\g_0>0,\mbox{ a.e. in }(0,\9)\times\rrd.
\end{equation}

In the following, we denote $\Phi(u(t,x))$ by $\Psi(t,x),$ $(t,x)\in(0,\9)\times\rrd.$

\begin{theorem}\label{t5} {\bf(Linearized uniqueness)} Under hypotheses {\rm(i)--(ii)}, for each $v_0\in L^1\cap L^\9$ and $T>0$, equation \eqref{e3.1} has at most one distributional solution $v\in L^{1}([0,T];L^1)\cap L^\9((0,T)\times\rrd).$  
\end{theorem}

\n{\bf Proof.} We shall proceed as in the proof of Theorem \ref{t1}. 
Namely, we set $v_1-v_2=v$ for  two solutions $v_1,v_2$ of \eqref{e3.2} 
By a similar argument as in the proof of Theorem \ref{t1}, we may assume that $v$ satisfies \eqref{add-e6}.
Then we get
\begin{equation}
\label{e3.3}
\barr{l}
v_t-\D(\Psi v)+{\rm div}(b(x,u)v)=0,\mbox{ a.e.  }t\in(0,T),\vsp 
v(0)=0.\earr
\end{equation}
For $y=\Gamma v$, we get
\begin{equation}
\label{e3.4}
\barr{l}
\dd\frac d{dt}\ y-\Gamma\D(\Psi v)+\Gamma\ {\rm div}(b(x,u)v)=0\vsp 
y(0)=0\earr
\end{equation}
and multiplying scalarly in $L^2$ with $v$, as above we get that
\begin{equation}
\label{e3.5}
\barr{l}
\dd\frac12\ \frac d{dt}\ |v(t)|^2_{-1}+\g_0|v(t)|^2_2
\le(|\Psi|_\9+|b|_\9)|v(t)|^{}_2|v(t)|^{}_{-1},\mbox{ a.e. }t\in(0,T).\earr
\end{equation}
This yields
$$\frac d{dt}\ |v(t)|^2_{-1}\le|v(t)|^2_{-1}\mbox{ a.e. }t\in(0,T),$$
and, therefore, $v\equiv0$, as claimed.

\section{Uniqueness in law of the McKean--Vlasov stochastic differential equations (SDEs)}
\setcounter{equation}{0}

Consider for $T\in(0,\9)$ and $u_0\in L^1\cap L^\9$ the McKean--Vlasov stochastic differential equation (SDE)
\begin{equation}\label{e4.1}
\barr{l}
dX(t)=b(X(t),u(t,X(t)))dt
+\dd\frac1{\sqrt{2}}\(\frac{\b(u(t,X(t)))}{u(t,X(t))}\)^{\frac12}dW(t),\vsp
\hfill 0\le t\le T,\\
u(0,\cdot)=\xi_0,\earr\end{equation}on $\rrd$. 
Here, $W(t),\ t\ge0,$ is an $(\calf_t)$-Brownian motion on a probability space $(\Omega,\calf,\mathbb{P})$ with normal filtration $\calf_t$, $t\ge0,$ $\xi_0:\Omega\to\rrd$ is $\calf_0$-measurable such that $$\mathbb{P}\circ\xi^{-1}_0(dx)=u_0(x)dx,$$ and $u(t,x)=\frac{d\mathcal{L}_{X(t)}}{dx}\,(x)$ is the Lebesgue density of the marginal law $\mathcal{L}_{X(t)}=\mathbb{P}\circ X(t)\1$ of the solution process $X(t)$, $t\ge0$. Here, a  solution process means an $(\calf_t)$-adapted process with $\pas$ continuous sample paths in  $\rrd$ solving \eqref{e4.1}.

\begin{theorem}\label{t4.1} Let $0<T<\9$ and let the above conditions {\rm(i)--(ii)} on  $b$ and $\b$ hold. Let $X(t)$, $t\ge0$, and $\wt X(t),$ $t\ge0$, be  two solutions to \eqref{e4.1} such that,~for 
$$u(t,\cdot):=\frac{d\mathcal{L}_{X(t)}}{dx},\ \ \wt u(t,\cdot):=\frac{d\mathcal{L}_{\wt X(t)}}{dx},$$we have
\begin{equation}\label{e4.2}
u,\wt u\in L^\9((0,T)\times\rrd).\end{equation}Then $X$ and $\wt X$ have the same laws, i.e., $\mathbb{P}\circ X\1=\mathbb{P}\circ\wt X^{-1}.$\end{theorem}

\n{\bf Proof.} By It\^o's formula, both $u$ and $\wt u$ satisfy the (nonlinear) \FP\ equation \eqref{e1.1} in the sense of Schwartz distributions. Hence, by Theorem \ref{t1}, $u=\wt u$. Furthermore, again by It\^o's formula, $\mathbb{P}\circ X\1$ and $\mathbb{P}\circ\wt X\1$ satisfy the martingale problem with the initial condition $u_0dx$ for the linear Komogorov operator
$$L_u:=\Phi(u)\D+b(\cdot,u)\cdot\nabla,$$where $\Phi(u)=\frac{\b(u)}u,$ $u\in\rr.$ Hence, by Theorem \ref{t5}, the assertion follows by Lemma 2.12 in \cite{4}. 

Here, for $s\in[0,T]$, the set $\mathcal{R}_{[s,T]}$, which appears in that lemma, is chosen to be the set of all narrowly continuous, probability measure-valued solutions of \eqref{e3.1} having for each $t\in[s,T]$ a density \mbox{$v(t,\cdot)\in L^\9$} with respect to  Lebesgue measure such that \mbox{$v\in L^\9((0,T)\times\rrd).$}\hfill$\Box$

\begin{remark}\label{e4.2}\rm We note that, by the narrow continuity, \eqref{e4.2} implies that, for every $t\in[0,T],$ $u(t,\cdot),\wt u(t,\cdot)\in L^\9$. This fact was used in the above proof.
\end{remark}

\n{\bf Acknowledgements.} This work was supported by the DFG through CRC 1283. We would like to thank Chengcheng Ling and Marco Rehmeier for pointing out a misprint that resulted in a gap in the proofs of two of the main results in an earlier version of this paper.


\begin{thebibliography}{nn}
	
	\bibitem{1} Barbu, V., R\"ockner, M., From nonlinear \FP\ equations to solutions of distribution dependent SDE, arXiv:1808.107062[math.PR].\vspace*{-2mm}  
	
	
	\bibitem{2az}  Barbu, V., R\"ockner, M., Probabilistic representation for solutions to nonlinear \FP\ equation, {\it SIAM J. Math. Anal.}, 50 (4) (2018), 4246-4260.\vspace*{-2mm}  
	
	\bibitem{2} Barbu, V., R\"ockner, M., The evolution to equilibrium of solutions to nonlinear \FP\ equations, arXiv:1904.082-91[math.PR].\vspace*{-2mm}     
	
	\bibitem{3}   Brezis, H., Crandall, M.G., Uniqueness of solutions of the initial-value problem for $u_t-\D\b(u)=0$, {\it J.~Math. Pures et Appl.}, 58 (1979), 153-163.\vspace*{-2mm}    	
	
	\bibitem{4} Trevisan, D., Well-posedness of multidimensional diffusion processes with weakly differentiable coefficients, {\it Electron. J. of Probab.}, Volume 21 (2016), paper no. 22, 41 pp.
	  
\end{thebibliography}
\end{document}